# SECTORIAL CONVERGENCE OF $U$-STATISTICS

By Anda Gadidov

*Kennesaw State University*

In this note we show that almost sure convergence to zero of symmetrized $U$-statistics indexed by a linear sector in $\mathbb{Z}_+^d$ is equivalent to convergence along the diagonal of $\mathbb{Z}_+^d$, as it is considered in Latała and Zinn [*Ann. Probab.* **28** (2000) 1908–1924]. Comparisons with similar results for sums of multi-indexed i.i.d. random variables are also made.

**1. Introduction.** Let $\mathbb{Z}_+^d, d \geq 1$, be the positive integer $d$-dimensional lattice points with coordinate-wise partial ordering $\leq$. A multi-index $(n_1, n_2, \ldots, n_d)$ in $\mathbb{Z}_+^d$ will be denoted $\mathbf{n}$. In particular, $n$ will be used to denote a $d$-tuple in which all indices are equal. For $0 < \theta < 1$, define the sector

$$\mathcal{S}_\theta^d = \left\{(i_1, i_2, \ldots, i_d) \in \mathbb{Z}_+^d : \theta < \frac{i_l}{i_k} < \frac{1}{\theta} \text{ for all } l, k = 1, \ldots, d\right\}.$$

Let $\{X_i\}$ be a sequence of i.i.d. random variables, $\{\varepsilon_i\}$ a Rademacher sequence independent of the $X_i$'s and $h$ a measurable function symmetric in its arguments. Denote $X_\mathbf{i} = (X_{i_1}, X_{i_2}, \ldots, X_{i_d})$ and $\varepsilon_\mathbf{i} = \varepsilon_{i_1} \varepsilon_{i_2} \cdots \varepsilon_{i_d}$.

Recently Latała and Zinn (2000) obtained necessary and sufficient conditions for the strong law of large numbers for symmetrized $U$-statistics with kernel $h$, $\gamma_n^{-1} \sum_{\mathbf{i} \in I_n} \varepsilon_\mathbf{i} h(X_\mathbf{i})$, where $I_n = \{\mathbf{i} \in \mathbb{Z}_+^d : \mathbf{i} \leq n, i_k \neq i_l, k \neq l\}$ and the normalizing sequence $\gamma_n$ satisfies some regularity conditions.

In this note we characterize the almost sure convergence to zero of $U$-statistics when the summation index set is the sector $\mathcal{S}_\theta^d$.

Convergence on rectangles with one vertex at the origin, in which the different indices go to infinity at their own pace (i.e., nonrestricted convergence), is studied in the context, of multi-sample $U$-statistics and necessary and suficient conditions have been obtained by McConnell (1987) for one and 2-sample $U$-statistics of order two.

Compared to the case considered by Latała and Zinn (2000), in which the limit is considered along the diagonal of $\mathbb{Z}_+^d$, the limiting index set is much richer in the case of a sector, but still more restrictive than the whole of $\mathbb{Z}_+^d$. Therefore, we might expect the results to be somewhere between McConnell's results and Latała and Zinn's.

We prove that sectorial convergence is equivalent to convergence along the diagonal of $\mathbb{Z}_+^d$. Almost sure convergence of the maxima of the normalized kernel is considered as well. The proofs are based on decoupling techniques, the Borel–Cantelli lemma and Lévy-type maximal inequalities. We also make use of some specific details in the proofs of Latała and Zinn (2000).

The main result is given in Section 2. Section 3 has some comments relating the result to similar treatments of sums of multi-indexed i.i.d. random variables.

Let us introduce some further notation. $X_{\mathbf{i}}^{dec} = (X_{i_1}^{(1)}, X_{i_2}^{(2)}, \ldots, X_{i_d}^{(d)})$, where $\{X_i^{(l)}\}$, $l=1,\ldots d$ are $d$ independent copies of $\{X_i\}$ and $\varepsilon_{\mathbf{i}}^{dec} = \varepsilon_{i_1}^{(1)} \varepsilon_{i_2}^{(2)} \cdots \varepsilon_{i_d}^{(d)}$, where $\varepsilon_i^{(l)}$, $l=1,\ldots,d$ are $d$ independent Rademacher sequences, independent of the $\{X_{\mathbf{i}}^{dec}\}$.

To avoid cluttered notation, $\max^\theta$ and $\sum^\theta$ will be used to indicate that the index set is restricted to the sector $\mathcal{S}_\theta^d$. Let us also convene that in all statements involving $h(X_{\mathbf{i}})$, the index $\mathbf{i}$ has all coordinates distinct, whereas in $h(X_{\mathbf{i}}^{dec})$ more than two coordinates can be identical.

**2. Sectorial convergence.** We assume, as in Latała and Zinn (2000), that the normalizing sequence $\gamma_{\mathbf{n}}$ satisfies the following conditions:

(2.1) $\quad \gamma_{\mathbf{n}}$ is nondecreasing in the coordinate-wise order on $\mathbb{Z}_+^d$,

(2.2) $\quad$ there exists $C>0$ such that $\gamma_{2n} \leq C\gamma_n$,

(2.3) $\quad \displaystyle\sum_{k \geq l} \frac{2^{dk}}{\gamma_{2^k}^2} \leq C \frac{2^{dl}}{\gamma_{2^l}^2} \quad$ for any $l=1,2,\ldots$.

We first consider the case of convergence in the sector of the normalized maxima. With no loss of generality, we will assume that $h$ is nonnegative.

THEOREM 2.1. *The following are equivalent:*

(2.4) $\quad \displaystyle\frac{1}{\gamma_{\mathbf{n}}} \max_{\mathbf{i} \leq \mathbf{n}}^\theta h(X_{\mathbf{i}}) \to 0 \quad$ *a.s. as $\mathbf{n} \to \infty$ in the sector $\mathcal{S}_\theta^d$,*

(2.5) $\quad \displaystyle\frac{1}{\gamma_n} \max_{\mathbf{i} \leq n} h(X_{\mathbf{i}}) \to 0 \quad$ *a.s.*

PROOF. Let $\{n_k\} \subset \mathbb{Z}_+$ such $1 < \liminf \frac{n_k}{n_{k-1}} \leq \limsup \frac{n_k}{n_{k-1}} < \infty$, and satisfying $n_{k-1} \leq \theta n_k < n_k$. One can define $n_k = p^k$, where $p > 1$ is an integer such that $\frac{1}{p} \leq \theta < \frac{1}{p-1}$.

Assume now that (2.4) holds. In particular,

$$\lim_{k\to\infty} \frac{1}{\gamma_{n_k}} \max_{\theta n_k < \mathbf{i} \le n_k}^{\theta} h(X_\mathbf{i}) = 0 \qquad \text{a.s.}$$

But $\{\mathbf{i} \in \mathbb{Z}_+^d : \theta n_k < \mathbf{i} \le n_k\} \subset \mathcal{S}_\theta^d$ and, therefore, by independence of the blocks, Borel–Cantelli lemma gives

$$\sum_k P\left\{\max_{\theta n_k < \mathbf{i} \le n_k} h(X_\mathbf{i}) > \epsilon \gamma_{n_k}\right\} < \infty$$

for all $\epsilon$, which, in turn, implies

$$\frac{1}{\gamma_{n_k}} \max_{\mathbf{i} \le (1-\theta)n_k} h(X_\mathbf{i}) \to 0 \qquad \text{a.s.}$$

Now (2.5) follows by the regularity of the sequence $\{\gamma_\mathbf{n}\}$.

Conversely, assume (2.5) holds and let $\mathbf{n} \in \mathcal{S}_\theta^d$. If $n^* = \max n_j$, we have

$$\frac{1}{\gamma_\mathbf{n}} \max_{\mathbf{i} \le \mathbf{n}}^{\theta} h(X_\mathbf{i}) \le \frac{1}{\gamma_\mathbf{n}} \max_{\mathbf{i} \le n^*}^{\theta} h(X_\mathbf{i}) \le \frac{C_1}{\gamma_{n^*}} \max_{\mathbf{i} \le n^*} h(X_\mathbf{i}) \to 0,$$

with the last inequality holding by the growth property (2.2). $\square$

A similar proof can show that the equivalence holds for the decoupled versions as well.

THEOREM 2.2. *The following are equivalent:*

(2.6) $\quad \dfrac{1}{\gamma_\mathbf{n}} \max_{\mathbf{i} \le \mathbf{n}}^{\theta} h(X_\mathbf{i}^{dec}) \to 0 \qquad$ *a.s. as* $\mathbf{n} \to \infty$ *in the sector* $\mathcal{S}_\theta^d$,

(2.7) $\quad \dfrac{1}{\gamma_n} \max_{\mathbf{i} \le n} h(X_\mathbf{i}^{dec}) \to 0 \qquad$ *a.s.*

Let us now look at the almost sure convergence of the $U$-statistics in the sector.

THEOREM 2.3. *The following are equivalent:*

(2.8) $\quad \dfrac{1}{\gamma_\mathbf{n}} \sum_{\mathbf{i} \le \mathbf{n}}^{\theta} \varepsilon_\mathbf{i} h(X_\mathbf{i}) \to 0 \qquad$ *a.s. as* $\mathbf{n} \to \infty$ *in the sector* $\mathcal{S}_\theta^d$,

(2.9) $\quad \dfrac{1}{\gamma_\mathbf{n}} \sum_{\mathbf{i} \le \mathbf{n}}^{\theta} \varepsilon_\mathbf{i}^{dec} h(X_\mathbf{i}^{dec}) \to 0 \qquad$ *a.s. as* $\mathbf{n} \to \infty$ *in the sector* $\mathcal{S}_\theta^d$,

(2.10) $\quad \dfrac{1}{\gamma_\mathbf{n}^2} \sum_{\mathbf{i} \le \mathbf{n}}^{\theta} h^2(X_\mathbf{i}) \to 0 \qquad$ *a.s. as* $\mathbf{n} \to \infty$ *in the sector* $\mathcal{S}_\theta^d$,

$$(2.11) \quad \frac{1}{\gamma_{\mathbf{n}}^2} \sum_{\mathbf{i} \leq \mathbf{n}}^{\theta} h^2(X_{\mathbf{i}}^{dec}) \to 0 \qquad \text{a.s. as } \mathbf{n} \to \infty \text{ in the sector } \mathcal{S}_{\theta}^d,$$

$$(2.12) \quad \frac{1}{\gamma_n} \sum_{\mathbf{i} \leq n} \varepsilon_{\mathbf{i}} h(X_{\mathbf{i}}) \to 0 \qquad \text{a.s.}$$

PROOF. $(2.8) \Rightarrow (2.10)$ and $(2.9) \Rightarrow (2.11)$ can be proved as in Cuzick, Giné and Zinn [(1995), Proposition 4.7].

$(2.10) \Rightarrow (2.12)$. We will actually show that (2.10) implies the convergence of the normalized sum of squares, which, by Latała and Zinn [(2000), Theorem 2], is equivalent to (2.12). If (2.10) holds, in particular,

$$\frac{1}{\gamma_{n_k}^2} \sum_{\theta n_k < \mathbf{i} \leq n_k}^{\theta} h^2(X_{\mathbf{i}}) \to 0 \qquad \text{a.s.}$$

But $\sum_{\theta n_k < \mathbf{i} \leq n_k}^{\theta} h^2(X_{\mathbf{i}}) = \sum_{\theta n_k < \mathbf{i} \leq n_k} h^2(X_{\mathbf{i}})$, and by independence of the blocks and the Borel–Cantelli lemma,

$$\sum_{k=1}^{\infty} P\left(\sum_{\theta n_k < \mathbf{i} \leq n_k} h^2(X_{\mathbf{i}}) > \epsilon \gamma_{n_k}^2\right) < \infty.$$

By the regularity of the $\gamma_{\mathbf{n}}$,

$$\frac{1}{\gamma_n} \sum_{\mathbf{i} \leq n} h^2(X_{\mathbf{i}}) \to 0 \qquad \text{a.s.,}$$

and (2.12) follows.

$(2.11) \Rightarrow (2.12)$. As before, it can be shown that (2.11) implies the a.s. convergence to zero of $\frac{1}{\gamma_{\mathbf{n}}^2} \sum_{\mathbf{i} \leq n} h^2(X_{\mathbf{i}}^{dec})$, which is equivalent to (2.12) by Latała and Zinn (2000), Theorem 2.

$(2.12) \Rightarrow (2.8)$ and $(2.12) \Rightarrow (2.9)$. In order to prove (2.8), it will be sufficient to show that

$$(2.13) \quad \sum_{k \geq 1} P\left\{\max_{\mathbf{n} \leq 2^k}^{\theta} \left|\sum_{\mathbf{i} \leq \mathbf{n}}^{\theta} \varepsilon_{\mathbf{i}} h(X_{\mathbf{i}})\right| > \epsilon \gamma_{2^k}\right\} < \infty.$$

Consider the $l^{\infty}$-space of vectors whose components are all possible sums $\sum_{\mathbf{i} \leq \mathbf{n}}^{\theta} \varepsilon_{\mathbf{i}} h(X_{\mathbf{i}})$, $\mathbf{n} \leq 2^k, \mathbf{n} \in \mathcal{S}_{\theta}^d$. By applying the decoupling inequality [Theorem 1 de la Peña and Montgomery-Smith (1995)] conditionally with respect to the Rademacher r.v. and $\{X_i\}$, respectively, there exists a constant $c_d > 0$, depending on $d$ only, such that

$$P\left\{\max_{\mathbf{n} \leq 2^k}^{\theta} \left|\sum_{\mathbf{i} \leq \mathbf{n}}^{\theta} \varepsilon_{\mathbf{i}} h(X_{\mathbf{i}})\right| > t\right\} = P\left\{\left\|\sum_{\mathbf{i} \leq \mathbf{n}}^{\theta} \varepsilon_{\mathbf{i}} h(X_{\mathbf{i}})\right\| > t\right\}$$

$$= \mathbb{E}_X P_\varepsilon \left\{ \left\| \sum_{\mathbf{i} \leq \mathbf{n}}^\theta \varepsilon_\mathbf{i} h(X_\mathbf{i}) \right\| > t \right\}$$

$$\leq \mathbb{E}_X c_d P_\varepsilon \left\{ \left\| \sum_{\mathbf{i} \leq \mathbf{n}, i_j \neq i_l}^\theta \varepsilon_\mathbf{i}^{dec} h(X_\mathbf{i}) \right\| > \frac{t}{c_d} \right\}$$

$$= \mathbb{E}_\varepsilon c_d P_X \left\{ \left\| \sum_{\mathbf{i} \leq \mathbf{n}}^\theta \varepsilon_\mathbf{i}^{dec} h(X_\mathbf{i}) \right\| > \frac{t}{c_d} \right\}$$

$$\leq c_d^2 P \left\{ \max_{\mathbf{n} \leq 2^k}^\theta \left| \sum_{\mathbf{i} \leq \mathbf{n}, i_j \neq i_l}^\theta \varepsilon_\mathbf{i}^{dec} h(X_\mathbf{i}^{dec}) \right| > \frac{t}{c_d^2} \right\}.$$

Therefore, (2.13) will follow if we show that

(2.14)
$$\sum_{k \geq 1} P \left\{ \max_{\mathbf{n} \leq 2^k}^\theta \left| \sum_{\mathbf{i} \leq \mathbf{n}}^\theta \varepsilon_\mathbf{i}^{dec} h(X_\mathbf{i}^{dec}) \right| > \epsilon \gamma_{2^k} \right\}$$
$$\leq 2^d \sum_{k \geq 1} P \left\{ \left| \sum_{\mathbf{i} \leq 2^k}^\theta \varepsilon_\mathbf{i}^{dec} h(X_\mathbf{i}^{dec}) \right| > \epsilon \gamma_{2^k} \right\} < \infty.$$

The above Lévy-type inequality can be proved iterately by applying Lévy's maximal inequality conditionally, $d$ times, as follows. For $1 \leq l \leq d$ and $1 \leq n_1, n_2, \ldots, n_d \leq 2^k$, let $P_l$ denote conditional probability given $\{\varepsilon_i^{(j)}, X_i^{(j)}, j \neq l\}$. Define the index sets $I_{l,k} = \{(\mathbf{i} \in \mathcal{S}_\theta^d : i_r \leq n_r, r \leq l, i_r \leq 2^k, r > l\}$ and $J_{l,k} = \{\mathbf{i} = (i_1, \ldots, i_{l-1}, i_{l+1}, \ldots, i_d) : i_r \leq n_r, r < l, i_r \leq 2^k, r > l\}$. For $1 \leq i \leq 2^k$, let $Y_i$ be a vector whose components are all possible sums,

$$\sum_{\mathbf{i} \in J_{l,k}}^\theta \varepsilon_{i_1}^{(1)} \cdots \varepsilon_{i_{l-1}}^{(l-1)} \varepsilon_{i_{l+1}}^{(l+1)} \cdots \varepsilon_{i_d}^{(d)} h(X_{i_1}^{(1)}, \ldots, X_i^{(l)}, \ldots, X_{i_d}^{(d)}).$$

Conditionally on $\{\varepsilon_i^{(j)}, X_i^{(j)}\}, j \neq l$, $\sum_{i \leq n_l} \varepsilon_i^{(l)} Y_i$ is a sum of independent and symmetric random vectors. Then, Lévy's maximal inequality gives

$$P \left\{ \max_{\mathbf{n} \leq 2^k}^\theta \left| \sum_{\mathbf{i} \in I_{l,k}}^\theta \varepsilon_\mathbf{i}^{dec} h(X_\mathbf{i}^{dec}) \right| > t \right\} = \mathbb{E} P_l \left\{ \max_{n_l \leq 2^k} \left\| \sum_{i \leq n_l} \varepsilon_i^{(l)} Y_i \right\| > t \right\}$$

$$\leq 2 \mathbb{E} P_l \left\{ \left\| \sum_{i \leq 2^k} \varepsilon_i^{(l)} Y_i \right\| > t \right\}$$

$$= 2 P \left\{ \max_{\mathbf{n} \leq 2^k}^\theta \left| \sum_{\mathbf{i} \in I_{l-1,k}}^\theta \varepsilon_\mathbf{i}^{dec} h(X_\mathbf{i}^{dec}) \right| > t \right\}.$$

Let us now define the sets $A_{k,1} = \{x \in E^d : h^2(x) \leq \gamma_{2^k}\}$ and for $l = 1, \ldots, d-1$, $A_{k,l+1} = \{x \in A_{k,l} : 2^{kl} \mathbb{E}_I h^2 I_{A_{k,l}}(x) \leq \gamma_{2^k}$ for all $I \subset \{1, 2, \ldots, d\}, \text{Card}(I) =$

$l\}$, as in Latała and Zinn (2000), Theorem 2. If (2.12) holds, then, by Theorem 2,

$$\sum_{k\geq 1} P\{\exists\, \mathbf{i} \leq 2^k : X_{\mathbf{i}}^{dec} \notin A_{k,d}\} < \infty.$$

Therefore, it only remains to prove

$$\sum_{k\geq 1} P\left\{\left|\sum_{\mathbf{i}\leq 2^k}^{\theta} \varepsilon_{\mathbf{i}}^{dec} h(X_{\mathbf{i}}^{dec}) I_{A_{k,d}}\right| > \epsilon \gamma_{2^k}\right\} \leq \sum_{k\geq 1} \frac{1}{\epsilon^2 \gamma_{2^k}^2} \mathbb{E}\left(\left|\sum_{\mathbf{i}\leq 2^k}^{\theta} \varepsilon_{\mathbf{i}}^{dec} h(X_{\mathbf{i}}^{dec}) I_{A_{k,d}}\right|^2\right)$$

$$\leq \sum_{k\geq 1} \frac{2^{kd}}{\gamma_{2^k}^2} \mathbb{E} h^2(X^{dec}) I_{A_{k,d}}(X^{dec}) < \infty,$$

the last inequality holding in view of (26) Latała and Zinn (2000).

Now (2.9) follows from (2.14) as well, and the conclusion follows. □

**3. Remarks and conclusions.** Questions regarding restricted convergence in the strong law of large numbers for sums of multi-indexed i.i.d. random variables arises quite naturally from the theory of convergence of multiple Fourier series or differentiability of multiple integrals. The problem of almost sure convergence when the index set is a partially ordered set in $\mathbb{Z}_+^d$ has been considered earlier by Smythe (1973, 1974), Gut (1983) and Klesov and Rychlik (1999) for sums of i.i.d. multi-indexed random variables.

Let $\mathcal{A}$ be a partially ordered subset of $\mathbb{Z}_+^d$. For $\alpha \in \mathcal{A}$, define $|\alpha| = \text{Card}\{\beta \in \mathcal{A} : \beta \leq \alpha\}$, and let

$$M(x) = \sum_{j\leq x} \text{Card}\{\alpha \in \mathcal{A} : |\alpha| = j\}, \qquad x \geq 0.$$

Smythe (1973) proved that, for a certain class of partially ordered sets, the Kolmogorov strong law of large numbers $|\beta|^{-1}\sum_{\beta\leq\alpha} X_\beta \to 0$ a.s. holds, if and only if $E(M(|X|)) < \infty$. In particular, for $\mathcal{A} = \mathbb{Z}_+^d$, $M(n) \sim n(\log n)^{d-1}$, and for $\mathcal{A} = \mathcal{S}_\theta^d$, $M(n) \sim n$.

Gut (1983) extended Smythe's result, proving that the Marcinkiewicz strong law of large numbers $|\mathbf{n}|^{-1/p}\sum_{\mathbf{i}\leq\mathbf{n}} X_{\mathbf{i}} \to 0$ a.s., $0 < p < 2$ and $EX = 0$ if $1 \leq p < 2$, holds in the sector $\mathcal{S}_\theta^d$ if and only if $E|X|^p < \infty$, which is exactly the necessary and sufficient condition for the classical strong law of large numbers. Klesov and Rychlik (1999) considered the case of almost sure convergence of normalized sums of i.i.d. random variables in a sector of $\mathbb{Z}_+^2$ with nonlinear boundaries. It turns out that in the case of sums of multi-indexed i.i.d. random variables, the strong law of large numbers is intrinsically related to the size of the index set.

Strong laws of large numbers for $U$-statistics are more complex since the summands display a nontrivial pattern of dependence. This is why a moment condition on the kernel of the $U$-statistic provides only a sufficient condition for the strong law of large numbers. However, the size of the limiting index set distinguishes between the various results. The necessary and sufficient conditions obtained by McConnell (1987) for $U$-statistics indexed by $\mathbb{Z}_+^2$ differ from the ones obtained by Latała and Zinn (2000). Notice that in Latała and Zinn (2000) the limiting index set is $\mathbb{Z}_+$. Moreover, the equivalence proved in this note supports this conclusion, since $M(n) \sim n$ for the sector $\mathcal{S}_\theta^d$, as well as for $\mathbb{Z}_+$.

**Acknowledgments.** The author wishes to thank Joel Zinn for some useful suggestions. Thanks are also due to the referee and the Editor for helpful comments.

Department of Mathematics
Kennesaw State University
1000 Chastain Road #1204
Kennesaw, Georgia 30144-5591
USA
e-mail: agadidov@kennesaw.edu
url: http://math.kennesaw.edu/~agadidov